\documentclass[12pt]{amsart}
\usepackage{amsmath, amssymb}

\author{}
\title{}
\date{}

\newcommand{\bbC}{{\mathbb C}}
\newcommand{\bbK}{{\mathbb K}}
\newcommand{\bbQ}{{\mathbb Q}}
\newcommand{\bbR}{{\mathbb R}}
\newcommand{\bbZ}{{\mathbb Z}}

\newcommand{\cA}{{\mathcal A}}
\newcommand{\cN}{{\mathcal N}}

\newcommand{\al}{\alpha}
\newcommand{\bbeta}{\mbox{\boldmath $\beta$}}
\newcommand{\br}{{\bf r}}

\newtheorem{corollary}{Corollary}
\newtheorem{lemma}{Lemma}

\begin{document}
\baselineskip=17pt

\vspace{10mm} 

\begin{center}
{\bf An effective lower bound for the height of algebraic numbers} \\
\vspace{5 mm}
by \\ 
\vspace{5 mm}
{\sc Paul M Voutier} (London) 
\end{center}

\vspace{7 mm}

{\bf 1. Introduction.}

Let $\al$ be a non-zero algebraic number of degree $d$ with 
\begin{displaymath}
f(X) = a_{d} \prod_{i=1}^{d} \left( X-\al_{i} \right)
\end{displaymath}
as its minimal polynomial over $\bbZ$ and $a_{d}$ positive. 
We shall define the Mahler measure of $\al$, $M(\al)$, by  
\begin{displaymath}
M(\al) = a_{d} \prod_{i=1}^{d} \max \left( 1, |\al_{i}| \right)
\end{displaymath}
and the absolute logarithmic height of $\al$, $h(\al)$, via 
the relationship 
\begin{displaymath}
h(\al) = \frac{\log M(\al)}{d}. 
\end{displaymath}

In 1933, D. H. Lehmer \cite{Lehmer} asked whether it is true that 
for every positive $\epsilon$ there exists an algebraic number  
$\al$ for which $1 < M(\al) < 1+\epsilon$. This question has since 
been reformulated as whether there exists a positive absolute 
constant $c_{0}$ such that $M(\al) > 1+c_{0}$; in this form the  
question is known as {\it Lehmer's problem}. The first progress 
was due to Schinzel and Zassenhaus \cite{SZ}, who proved in 1965 
that if $\al$ is not a root of unity then $\overline{|\al|} 
> 1+4^{-s-2}$ where $2s$ is the number of complex conjugates  
of $\al$. This implies that $M(\al) > 1+c_{1}/2^{d}$ for a 
positive absolute constant $c_{1}$. In 1971, Blanksby and Montgomery 
\cite{BM} used Fourier analysis to make a considerable refinement 
upon this first result. They proved that $M(\al) > 1+1/(52d\log(6d))$. 
In 1978, C.L. Stewart \cite{Stew} introduced a method from transcendental 
number theory to prove that $M(\al) > 1+1/(10^{4}d \log (d))$. 
While this result is a little weaker than the one due to Blanksby 
and Montgomery, the method used has since become quite important 
as it has produced the best results yet known, results which are 
significantly better than those previously known and bring us 
quite close to the conjectured lower bound.  

Dobrowolski \cite{Dob} was able to use Stewart's method to show that 
for each $\epsilon > 0$, there exists a positive integer $d(\epsilon)$ 
such that for all $d \geq d(\epsilon)$, 
\begin{displaymath}
M(\al) > 1 + (1-\epsilon) { \left( \frac{\log \log d}{\log d} \right) }^{3}.  
\end{displaymath}

Moreover, his proof can be made effective. Dobrowolski claims that,  
for all $d \geq 3$, 
\begin{displaymath}
M(\al) > 1 + \frac{1}{1200} \left( \frac{\log \log d}{\log d} \right)^{3}.  
\end{displaymath}

Dobrowolski achieves his results by constructing a polynomial,  
$F(X) \in \bbZ[X]$, of small height which is divisible by $f(X)^{T}$,   
for a certain integer $T$ which depends on $d$, and considering the 
norm of $F(\al^{p})$ at certain primes $p$. The big improvement of 
his lower bound over previous ones arises by a clever use of Fermat's 
little theorem to replace the trivial lower bound 
$1 \leq |\cN_{\bbQ(\al)/\bbQ}(F(\al^{p}))|$ by 
$p^{dT} \leq |\cN_{\bbQ(\al)/\bbQ}(F(\al^{p}))|$.  

In 1981, Cantor and Straus \cite{CS} showed that instead of using 
the auxiliary function $F(X)$, one could develop a version of 
Dobrowolski's proof by considering the determinant of a certain 
matrix. Fermat's little theorem still permits a non-trivial lower 
bound for the absolute value of this determinant while Hadamard's 
inequality gives an upper bound which involves $M(\al)$. Comparing 
these bounds, as in Dobrowolski's proof, yields the desired lower 
bound upon judicious choice of certain parameters. Their proof had 
the advantages of being simpler than Dobrowolski's as well as yielding 
an improved result: they were able to replace the constant $1-\epsilon$ 
by $2-\epsilon$. Louboutin \cite{Loub}, in 1983, used a more refined 
selection of parameters to improve this constant to $9/4-\epsilon$. 
In 1988, Meyer \cite{Meyer} showed that one can also obtain this result 
from Dobrowolski's method using an auxiliary function.  

In addition to these results there is also a very important result 
due to Smyth \cite{Smyth}. In 1971, he proved that if $\al^{-1}$ is 
not a conjugate of $\al$ then $M(\al) \geq \al_{1} = 1.32471 \ldots$. 
This number, $\al_{1}$, is the real root of the polynomial $X^{3}-X-1$ 
and also the smallest Pisot number. One consequence of his result is 
the positive solution of Lehmer's problem when $d$ is odd with 
$c_{0}=0.32471 \ldots$, since reciprocal polynomials (those with 
$\al^{-1}$ as a root whenever $\al$ is a root) of odd degree 
have -1 as a root and are hence reducible. Therefore, it is now only 
necessary to consider even values of $d$.  

Finally, let us mention some computational work done on this 
problem. Lehmer performed considerable computations in the early 
thirties. The smallest value of $M(\al)$ he was able to find arose 
from the example $\al=M(\al)=1.1762808 \ldots$ which is a root of 
the polynomial 
\begin{displaymath}
X^{10}+X^{9}-X^{7}-X^{6}-X^{5}-X^{4}-X^{3}+X+1.   
\end{displaymath}

This value of $M(\al)$ is still the smallest known --- indeed, 
it is widely believed that this is the minimum value of $M(\al)$. 
Interestingly, this $\al$ is also the smallest known Salem number.

Boyd \cite{Boyd1,Boyd2} has conducted extensive calculations on 
this problem and determined all $\al$ of degree at most 20 with 
$M(\al) \leq 1.3$ without finding smaller values. Thus, we need 
only consider even values of $d \geq 22$.

We mentioned before that Dobrowolski's proof can be made effective. 
So can those of Cantor and Straus and Louboutin. In this article,  
we produce such an effective result using the method of Cantor and 
Straus. Using their proof has the advantage of simplicity over the 
other two methods. Moreover, although it yields a result which is 
asymptotically weaker than Louboutin's lower bound, for small $d$ 
the methods give rise to the same choice of parameters and so we do 
not lose anything by this choice. Let us now state our result. 

\vspace{3.0mm}

{\bf Theorem.} {\it Suppose $\al$ is a non-zero algebraic number 
of degree $d$ which is not a root of unity. If $d \geq 2$ then} 
\begin{displaymath}
dh(\al) = \log M(\al) > \frac{1}{4} 
			\left( \frac{\log \log d}{\log d} \right)^{3}.  
\end{displaymath}

\vspace{3.0mm}

There is a conjecture related to Lehmer's problem due to Schinzel 
and Zassenhaus. They ask whether there exists a positive absolute 
constant $c_{2}$ such that the maximum of the absolute values of 
the conjugates of $\al$, denoted $\overline{|\al|}$, always 
satisfies $\overline{|\al|} \geq 1+c_{2}/d$ when $\al$ is a 
non-zero algebraic number which is not a root of unity. Notice 
that $1+\log(M(\al)) < M(\al) \leq \overline{|\al|}^{d}$ so 
$1+\log(M(\al))/d < \overline{|\al|}$. Therefore, Smyth's result 
implies that this conjecture is true when $\al^{-1}$ is not an  
algebraic conjugate of $\al$. When $\al^{-1}$ and $\al$ are 
algebraic conjugates, we can replace the inequality above by 
$1+\log(M(\al)) < \overline{|\al|}^{d/2}$ for such $\al$, 
and so our theorem yields the following corollary. 

\begin{corollary}
\label{cor:mycor1}
Let $\al$ be as in the theorem above. Then 
\begin{displaymath}
\overline{|\al|} 
> 1 + \frac{1}{2d} \left( \frac{\log \log d}{\log d} \right)^{3}.  
\end{displaymath}
\end{corollary}

The best known result on this problem is that 
\begin{displaymath}
\overline{|\al|} 
> 1 + \left( \frac{64}{\pi^{2}} - \epsilon \right) 
      \frac{1}{d} \left( \frac{\log \log d}{\log d} \right)^{3},   
\end{displaymath}
for $d \geq d(\epsilon)$. This result was proven in 1993 by 
Dubickas \cite{Dub}. 

For $2 \leq d \leq 2300$, there is a better result than 
Corollary~\ref{cor:mycor1} which is due to Matveev \cite{Matveev}: 
for $d \geq 1$, 
\begin{displaymath}
\overline{|\alpha|} > \exp \left( \frac{3\log (d/2)}{d^{2}} \right). 
\end{displaymath}

In applications, it will often be more convenient to use a simpler 
form of these results. Therefore, we record the following corollary 
which follows from the theorem, the work of Boyd \cite{Boyd1,Boyd2} 
and Smyth \cite{Smyth} and the inequality $1+\log(M(\al)) 
< {\overline{|\al|}}^{d/2}$ when $\al$ and $\al^{-1}$ are algebraic 
conjugates along with the fact that, as shall be demonstrated in the 
course of the proof of the theorem, we may replace $1/4$ by $0.56$ for 
$22 \leq d \leq 190$. Let us note that, at the expense of more effort 
and further complications, we could prove such a result for $d > 190$ 
too. 

\begin{corollary}
\label{cor:mycor2}
Suppose $\al$ is a non-zero algebraic number of degree $d$ which 
is not a root of unity. If $d \geq 2$ then 
\begin{displaymath}
dh(\al) = \log M(\al) > \frac{2}{(\log (3d))^{3}} 
\end{displaymath}
and 
\begin{displaymath}
\overline{|\al|} > 1 + \frac{4}{d(\log (3d))^{3}} 
\end{displaymath}
\end{corollary}

\vspace{3 mm}

{\bf 2. Preliminary Lemmas.} 

Let us begin with some lemmas. 

\begin{lemma}
\label{lem:prod}
If $K$ is a positive integer then 
\begin{displaymath}
\prod_{k=0}^{K-1} (k!) 
\geq \exp \left( \frac{K^{2} \log K}{2} - \frac{3K^{2}}{4} \right). 
\end{displaymath}
\end{lemma}

\begin{proof}
For $K=1$, the left-hand side is $1$, whereas the right-hand 
side is $0.47\ldots$. So we may assume that $K \geq 2$. 

In the course of proving their Lemme 8, Laurent, Mignotte and 
Nesterenko \cite{LMN} showed that 
\begin{displaymath}
\prod_{k=0}^{K-1} (k!) 
\geq \exp \left( \frac{(K^{2}-K) \log (K-1)}{2} - \frac{3(K^{2}-K)}{4} 
+ \frac{K \log (2\pi (K-1)/\sqrt{e})}{2} - \frac{\log K}{12} \right). 
\end{displaymath}

By means of the relation $\log (K-1) = \log (K) + \log (1-1/K)$ 
and the series expansion for $\log (1-x)$, we can replace the  
two terms of the form $\log (K-1)$ by $\log (K) - 7/(5K)$ for 
$K \geq 2$ to obtain 
\begin{displaymath}
\prod_{k=0}^{K-1} (k!) 
\geq \exp \left( \frac{K^{2} \log K}{2} - \frac{3K^{2}}{4}  
		 + K\log (\sqrt{2\pi}e^{-1/5}) - \frac{\log K}{12} \right). 
\end{displaymath}

The lemma now follows. 
\end{proof}

\begin{lemma}
\label{lem:primesum}
Let $p_{1}=2, p_{2}=3, \ldots$ denote the prime numbers in increasing 
order. 

\noindent
{\rm (i)} For $S \geq 13$, 
\begin{displaymath}
\theta(p_{S}) = \sum_{i=1}^{S} \log p_{i} \geq S \log S.   
\end{displaymath}

\noindent
{\rm (ii)} For $i \geq 20$, 
\begin{displaymath}
p_{i} < i(\log i + \log \log i - 1/2). 
\end{displaymath}

\noindent
{\rm (iii)} For $S \geq 9$, 
\begin{displaymath}
\sum_{i=1}^{S} p_{i} \leq 0.564S^{2} \log S.  
\end{displaymath}
\end{lemma}

\begin{proof}
(i) This is Th\'{e}or\`{e}me 4 of \cite{Robin}. 

(ii) This is equation (3.11) in the statement of Theorem~3 
of \cite{RS}. We will use it to prove part (iii). 

(iii) One can easily check that this is true for 
$9 \leq S \leq 19$. In the course of this calculation we 
find that $p_{1} + \cdots + p_{19}=568$, a fact that we 
shall now use to prove the inequality in general. 

For $S \geq 20$, we have
\begin{eqnarray*}
\sum_{i=1}^{S} p_{i} 
&  <   & 568 + \sum_{i=20}^{S} ( i \log i + i \log \log i - i/2 ) \\ 
& \leq & 568 
	 + \int_{i=20}^{S} \left( i \log i + i \log \log i - i/2 \right)  di 
	 + S \log S + S \log \log S - S/2.  
\end{eqnarray*}

Knowing that 
\begin{displaymath}
\int x \log x + x \log \log x - x/2 \, dx 
= \frac{x^{2}}{2} \left( \log x + \log \log x - 1 \right) - li(x^{2})/2,   
\end{displaymath}
that $li(x^{2}) > x^{2}/(2 \log (x))$ and that $li(400) = 85.417 \ldots$, 
we find that 
\begin{eqnarray*}
\sum_{i=1}^{S} p_{i} 
< \frac{S^{2}}{2} \left( \log S + \log \log S - 1 - \frac{1}{2 \log S} \right) 
  + S \log S + S \log \log S - S/2 - 7.8.  
\end{eqnarray*}

To prove the desired inequality we want to consider 
\begin{displaymath}
f(S) = \frac{\log \log S - 1 - 1/(2 \log S)}{\log S} 
       +  \frac{2(S \log S + S \log \log S - S/2 - 7.8)}{S^{2} \log S}. 
\end{displaymath}

Taking the derivative of $f(S)$, we find that its maximum for 
$S \geq 20$, which occurs at $S=2803.26\ldots$, is less than 0.128.  
This implies our result. 
\end{proof}

\vspace{3 mm}

We now want to give a technical lemma for later use. 

\begin{lemma}
\label{lem:positive}
Suppose $d > 10000$, $1.2 \leq k_{1} \leq 1.6$ and $1.2 \leq s_{1}$, 
then 
\begin{displaymath}
\left( k_{1}^{2} + 2s_{1} \right) \frac{\log \log \log d}{\log \log d}
+ \frac{\left( 0.89 + \log k_{1} \right) k_{1}^{2}+2.39s_{1}}{\log \log d}  
- (k_{1}^{2} + 2s_{1}) 
> - \frac{2k_{1}s_{1} \log s_{1} \log d}{(\log \log d)^{2}}.  
\end{displaymath}
\end{lemma}

\begin{proof}
Notice that $0.89 + \log 1.2 > 1.07$ and put 
\begin{displaymath}
f_{1}(d)= 1 - \frac{\log \log \log d}{\log \log d} 
	  - \frac{1.07}{\log \log d} \mbox{ and } 
f_{2}(d)= 2 - \frac{2 \log \log \log d}{\log \log d} 
	  - \frac{2.39}{\log \log d}. 
\end{displaymath}          

We want show that 
\begin{displaymath}
g(d,k_{1},s_{1}) = \frac{2k_{1} s_{1} \log s_{1} \log d}{(\log \log d)^{2}} 
		       - f_{1}(d)k_{1}^{2} - f_{2}(d)s_{1} 
\end{displaymath}                   
is positive for $d > 10000, 1.2 \leq k_{1} \leq 1.6$ 
and $1.2 \leq s_{1}$. 

We start by noting that 
\begin{displaymath}
f_{1}'(d) = \frac{7 + 100 \log \log \log d}{100 d \log d (\log \log d)^{2}} 
\hspace*{5 mm} \mbox{ and } \hspace*{5 mm}
f_{2}'(d) = \frac{39 + 200 \log \log \log d}{100d \log d (\log \log d)^{2}},  
\end{displaymath}          
so $f_{1}$ and $f_{2}$ are increasing functions for 
$d \geq \exp \exp 1 = 15.15 \ldots$. We shall also need 
to know that $(\log d)/(\log \log d)^{2}$ is increasing for 
$d \geq \exp \exp 2 = 1618.17 \ldots$.

We will prove that $g(d,k_{1},s_{1})$ is positive in the 
desired domain by considering six different ranges of $d$. 

If $10000 < d \leq 10^{8}$ then $(\log d)/(\log \log d)^{2} > 1.865$, 
$f_{1}(d) < 0.27$ and $f_{2}(d) < 0.45$ so 
$g(d,k_{1},s_{1}) > g_{1}(k_{1},s_{1}) 
= 3.73 k_{1}s_{1} \log s_{1} - 0.27 k_{1}^{2} - 0.45 s_{1}$ 
for these $d, k_{1} > 0$ and $s_{1} > 1$. Notice that 
$g_{1}(1.2,1.2) > 0.05$ and that $g_{1}(k_{1},1.2)$ is a 
quadratic polynomial which is positive between its two roots 
at $0.978\ldots$ and $2.044\ldots$. These two facts combined 
with the fact that $(\partial/ \partial s_{1}) g_{1}(k_{1},s_{1}) 
=3.73k_{1} \log s_{1} + 3.73 k_{1} - 0.45$ is positive for our 
range of $k_{1}$ and $s_{1}$ proves that $g_{1}$, and hence $g$, 
is positive for these values of $d$.

In the next four intervals we proceed similarly. 

For $10^{8} < d \leq 10^{12}$, we use 
$g_{1}=4.34k_{1}s_{1}\log s_{1} - 0.32k_{1}^{2} - 0.56s_{1}$ and 
find that $g_{1}(1.2,1.2) > 0.0066$, the roots of $g_{1}(k_{1},1.2)$ 
are $1.165\ldots$ and $1.801\ldots$ and 
$(\partial/ \partial s_{1}) g_{1}(k_{1},s_{1})$ is positive for 
our range of $k_{1}$ and $s_{1}$ to prove the lemma. 

For $10^{12} < d \leq 10^{18}$, we use 
$g_{1}=5.01k_{1}s_{1}\log s_{1} - 0.36k_{1}^{2} - 0.66s_{1}$; 
$g_{1}(1.2,1.2) > 0.0049$, the roots of $g_{1}(k_{1},1.2)$ 
are $1.179\ldots$ and $1.865\ldots$ and again 
$(\partial/ \partial s_{1}) g_{1}(k_{1},s_{1})$ is positive 
for our range of $k_{1}$ and $s_{1}$ which, as before establishes 
the lemma for this range of $d$. 

For $10^{18} < d \leq 10^{30}$, we use 
$g_{1}=5.97k_{1}s_{1}\log s_{1} - 0.41k_{1}^{2} - 0.76s_{1}$.  
The roots of $g_{1}(k_{1},1.2)$ are $1.033\ldots$ and $2.152\ldots$,  
$g_{1}(1.2,1.2) > 0.064$ and 
$(\partial/ \partial s_{1}) g_{1}(k_{1},s_{1})$ is positive 
for our range of $k_{1}$ and $s_{1}$. This proves the lemma for  
this range of $d$. 

For $10^{30} < d \leq 10^{100}$, we use 
$g_{1}=7.7k_{1}s_{1}\log s_{1} - 0.5k_{1}^{2} - 0.94s_{1}$ and 
find that $g_{1}(1.2,1.2) > 0.17$, the roots of $g_{1}(k_{1},1.2)$ 
are $0.921\ldots$ and $2.447\ldots$ and again \\ 
$(\partial/ \partial s_{1}) g_{1}(k_{1},s_{1})$ is positive for 
our range of $k_{1}$ and $s_{1}$ to prove the lemma. 

For $d > 10^{100}$, we have $f_{1}(d) < 1$, $f_{2}(d) < 2$ 
and $\log d/(\log \log d)^{2} > 7.78$. Therefore, 
$g(k_{1},s_{1},d) > h(k_{1},s_{1})
=15.56k_{1}s_{1} \log s_{1} -k_{1}^{2}-2s_{1}$ for $k_{1} > 0$,  
$s_{1} > 1$ and our range of $d$. As in the previous cases, we use 
the facts that $h(1.2,1.2)=0.245 \ldots$, the partial derivative 
of $h$ with respect to $s_{1}$ is positive for $k_{1} > 1$ and 
$s_{1} > 1$, and the roots of $h(k_{1},1.2)$ are $0.996 \ldots$ 
and $2.407 \ldots$ to prove the lemma in this remaining case. 
\end{proof}

\vspace{3.0mm}

To allow us to get a good constant in our theorem, we shall make 
use of the following result bounding the absolute value of the 
discriminant of a number field in terms of its degree. 

\begin{lemma}
\label{lem:disc}
Let $\bbK$ be an algebraic number field of degree $d$ with $D_{\bbK}$ 
as its discriminant. Then 
\begin{displaymath}
\log \left| D_{\bbK} \right| > 3.108d-8.6d^{1/3}.
\end{displaymath}
\end{lemma}

\begin{proof}
This is a result of Odlyzko, see equation (22) of \cite{Poitou}. 
\end{proof}

\vspace{3 mm}

The last lemma we give before introducing the ideas of Cantor and 
Straus will help us to prove that a certain determinant is not zero. 

\begin{lemma}
\label{lem:rootofunity}
Suppose $\al$ is a non-zero algebraic integer of degree $d$ 
with $\al=\al_{1}, \ldots, \al_{d}$ as its conjugates. If 
there exist positive rational numbers $r$ and $s$ such that 
$\al_{i}^{r} = \al_{j}^{s}$ then either $r=s$ or $\al$ is 
a root of unity. 
\end{lemma}

\begin{proof}
This is Lemma 2(i) of \cite{Dob}.
\end{proof}

\vspace{3 mm}

We now come to the work of Cantor and Straus. 

Let $v_{0}(\beta) = \left( 1,\beta, \beta^{2}, \ldots, \beta^{n-1} \right)^{t}$
and
\begin{displaymath}
v_{i}(\beta) = \frac{1}{i!} \frac{d^{i}}{d \beta^{i}} v_{0}(\beta) 
= \left( {0 \choose i} \beta^{-i}, {1 \choose i} \beta^{1-i}, \ldots, 
	 {n-1 \choose i} \beta^{n-1-i} \right)^{t} \mbox{ for $i \geq 1$,} 
\end{displaymath}
where we set ${h \choose i} = 0$ if $h < i$ and $h \in \bbZ$. 

Suppose that 
\begin{displaymath}
\bbeta = \left( \beta_{1}, \ldots, \beta_{m} \right) \in \bbC^{m} 
\mbox{ and } 
\br = \left( r_{1}, \ldots, r_{m} \right) \in \bbZ_{>0}^{m}.  
\end{displaymath}

We put $n= \sum_{j=1}^{m} r_{j}$ and define the confluent 
Vandermonde determinant $V(\bbeta,\br)$ to be the 
determinant of the $n \times n$ matrix whose columns are the 
vectors $v_{i} \left( \beta_{j} \right)$ where $0 \leq i \leq r_{j}-1$ and
$1 \leq j \leq m$. 

This $n \times n$ matrix bears some relation to the matrix 
associated with the system of linear equations from which 
Dobrowolski \cite{Dob} constructs his auxiliary function. 
In addition, a resemblance which this matrix bears to the 
Vandermonde matrix yields a particularly elegant formula 
for its determinant. 

\begin{lemma}
\label{lem:det}
Suppose that $\bbeta$ and $\br$ are as above. Then 
\begin{displaymath}
V(\bbeta, \br) 
= \pm \prod_{1 \leq i < j \leq m} { \left( \beta_{i} - \beta_{j} \right) 
				  }^{r_{i}r_{j}}. 
\end{displaymath}
\end{lemma}

\begin{proof}
This is Lemma 1 of \cite{CS}. 
\end{proof}

\vspace{3 mm}

Let us now see how to use this determinant to prove our theorem. 

Let $k$ and $s$ be two positive integers. Put $p_{0}=1$ and let 
$p_{1}, \ldots, p_{s}$ be the first $s$ prime numbers. We define 
\begin{displaymath}
\bbeta = \left( \alpha_{1}^{p_{0}}, \ldots, \alpha_{d}^{p_{0}},           
		\alpha_{1}^{p_{1}}, \ldots, \alpha_{d}^{p_{1}}, \ldots,    
		\alpha_{1}^{p_{s}}, \ldots, \alpha_{d}^{p_{s}} \right) 
\mbox{ and }                     
\br = \left( k,\ldots,k,1,1,\ldots,1 \right), 
\end{displaymath}
where the first $d$ components of $\br$ are $k$'s and the last 
$sd$ components are $1$'s. Notice that $m=(s+1)d$ and $n=d(k+s)$. 

From the previous three lemmas, we shall determine a lower 
bound for $|V(\bbeta,\br)|$. 

\begin{lemma}
\label{lem:lowbnd}
Suppose that $\al$ is a non-zero algebraic integer of degree 
$d \geq 2$ over $\bbQ$ which is not a root of unity and 
$\bbQ \left( \al^{p} \right) = \bbQ(\al)$ for all primes $p$. Moreover, suppose
that $\bbeta$ and $\br$ are as defined above. Then 
\begin{displaymath}
\left| V(\bbeta,\br) \right|^{2}  
\geq { \left( p_{1} \cdots p_{s} \right) }^{2dk} 
\exp \left( \left( 3.108 - 8.6d^{-2/3} \right) d(k^{2}+s) \right). 
\end{displaymath}     
\end{lemma}

\begin{proof}
Let us start by showing that $V=V(\bbeta,\br) \neq 0$ 
with our choice of $\bbeta$ and $\br$. From our expression 
for $V$ in Lemma~\ref{lem:det} and our definition of $\bbeta$, 
we see that $V=0$ if and only if there exist integers $i,j,k$ 
and $\ell$ such that $\al_{i}^{p_{k}}=\al_{j}^{p_{\ell}}$. Clearly  
$i \neq j$, for otherwise $\al$ is a root of unity (notice that 
due to the form of the expression in Lemma~\ref{lem:det} along with 
our choice of  $\bbeta$ and $\br$, if $i=j$ then $k \neq \ell$). 
Next, by Lemma~\ref{lem:rootofunity}, $k=\ell$ unless $\al$ is a root 
of unity. So we need only consider the case $\al_{i}^{p}=\al_{j}^{p}$ 
for some prime $p$. 

Define the polynomial $f_{p}(X)$ by 
\begin{displaymath}
f_{p}(X) = \prod_{j=1}^{d} \left( X - \al_{j}^{p} \right). 
\end{displaymath} 

By Lemme~7.1.1 of \cite{BP}, $f_{p}(X) \in \bbZ[X]$ and is either the minimal
polynomial of $\al^{p}$ or a power of 
this minimal polynomial. Now if $\al_{i}^{p}=\al_{j}^{p}$ 
for some pair of distinct integers $i$ and $j$, then 
$f_{p}(X)$ has multiple roots so it must be a power of the 
minimal polynomial of $\al^{p}$. But this implies that 
$[\bbQ(\al^{p}):\bbQ] < [\bbQ(\al):\bbQ]$, which contradicts 
one of our hypotheses. 

Thus, in what follows, we can suppose that $V$ is non-zero. 

With this in mind, let us now obtain a lower bound for $V^{2}$. 
From the expression in Lemma~\ref{lem:det} along with our 
definitions of $\bbeta$ and \br, we see that 
\begin{eqnarray*}
V^{2} & = & { \left( \prod_{1 \leq i < j \leq d} 
		     \left( \al_{i} - \al_{j} \right) \right) }^{2k^{2}}
	    { \left( \prod_{\ell=1}^{s} \prod_{1 \leq i,j \leq d} 
		     \left( \al_{i}^{p_{\ell}} - \al_{j} \right) \right) }^{2k}  
	      \left( \prod_{\ell=1}^{s} \prod_{1 \leq i < j \leq d} 
		     { \left( \al_{i}^{p_{\ell}} - \al_{j}^{p_{\ell}} \right) }^{2} 
	      \right) \\
      &   & { \left( \prod_{1 \leq \ell_{1} < \ell_{2} \leq s} 
		     \prod_{1 \leq i,j \leq d} 
		     \left( \al_{i}^{p_{\ell_{1}}} - \al_{j}^{p_{\ell_{2}}} 
		     \right) \right) }^{2}.  
\end{eqnarray*}

Let us denote these products by $A_{1}, A_{2}, A_{3}$ and $A_{4}$ 
according to their order of appearance. Notice that each of the 
$A_{i}$'s,  and hence $V^{2}$, is a symmetric function in the 
$\alpha_{i}$'s and thus a rational integer since the $\alpha_{i}$'s 
are algebraic integers. We shall now determine integers which divide 
these $A_{i}$'s.

We first consider $A_{2}$. Let $f(X)$ be the minimal polynomial over 
$\bbQ$ of $\alpha$ and notice, by Fermat's little theorem, that 
\begin{displaymath}
f \left( X^{p} \right) \equiv f(X)^{p} \bmod p, 
\end{displaymath}
for any prime $p$. Therefore 
\begin{displaymath}
f \left( \alpha_{i}^{p} \right) \equiv f \left( \alpha_{i} \right)^{p} \equiv 0 \bmod p, 
\end{displaymath}
for $1 \leq i \leq d$. 

Thus $p^{d}$ divides 
\begin{displaymath}
\prod_{i=1}^{d} f \left( \alpha_{i}^{p} \right)
= \prod_{i=1}^{d} \prod_{j=1}^{d} \left( \alpha_{i}^{p} - \alpha_{j} \right)  
\end{displaymath}
and therefore $(p_{1} \cdots p_{s})^{2dk}$ divides $A_{2}$. This will 
provide the main term in our lower bound for $V^{2}$.  

Now let us examine $A_{1}$ and $A_{3}$. $A_{1}$ is simply 
$\mbox{disc}(\al)^{k^{2}}$ and the inner product in $A_{3}$ 
is $\mbox{disc}(\al^{p_{\ell}})$. Since $\al$ is an algebraic 
integer and since we assumed that $\bbQ(\al)=\bbQ \left( \al^{p_{\ell}} \right)$, 
these discriminants are both divisible by $D_{\bbQ(\al)}$, the 
discriminant of $\bbQ(\al)$. Therefore, $D_{\bbQ(\al)}^{k^{2}+s}$ 
divides $A_{1}A_{3}$. 

We now consider
$$
A_{4} = { \left( \prod_{1 \leq \ell_{1} < \ell_{2} \leq s} 
		     \prod_{1 \leq i,j \leq d} 
		     \left( \al_{i}^{p_{\ell_{1}}} - \al_{j}^{p_{\ell_{2}}} 
		     \right) \right) }^{2}.
$$

The reader will notice that the product $A_{4}$ contains a 
large number of terms in comparison with the other three terms 
and thus should make a large contribution to a lower bound for 
$|V^{2}|$; however, we have been unable to determine a lower 
bound for $|A_{4}|$ other than the trivial one: $\left| A_{4} \right| \geq 1$.

Combining these results with our assumption that $V \neq 0$, 
we find that 
\begin{displaymath}
\left| V^{2} \right| 
\geq { \left( \prod_{i=1}^{s} p_{i} \right) }^{2dk}
     { \left| D_{\bbQ(\al)} \right| }^{k^{2}+s}. 
\end{displaymath}     

Applying Lemma~\ref{lem:disc} completes the proof of the lemma.  
\end{proof}

\vspace{3 mm}

We will also need an upper bound. For this, we shall use 
Hadamard's inequality.  

\begin{lemma}
\label{lem:upbnd}
Suppose that $\bbeta$ and $\br$ are as above. Then 
\begin{displaymath}
\left| V(\bbeta,\br) \right|^{2}  
\leq { \left( \frac{n^{k^{2}+s}}
		   {\displaystyle \prod_{i=0}^{k-1} (2i+1)(i!)^{2}} 
       \right) }^{d} M(\alpha)^{2n(k+p_{1}+\cdots+p_{s})}.          
\end{displaymath}     
\end{lemma}

\begin{proof}
By Hadamard's inequality, we have 
\begin{displaymath}
\left| V(\bbeta,\br) \right|^{2}  
\leq \prod_{j=1}^{m} \prod_{i=0}^{r_{j}-1} \left| v_{i} \left( \beta_{j} \right) \right|^{2} 
  =  \left( \prod_{j=1}^{d} \prod_{i=0}^{k-1} \left| v_{i} \left(\beta_{j} \right) \right|^{2} \right)
     \times \left( \prod_{j=d+1}^{m} \left| v_{0} \left( \beta_{j} \right) \right|^{2} \right).
\end{displaymath}         

By the definition of $v_{0}$, 
$\left| v_{0} \left( \beta_{j} \right) \right|^{2} \leq n \max \left( 1, \left| \beta_{j} \right|  \right)^{2(n-1)}$.
Using this along with our definition of $\bbeta$, we find that 
\begin{eqnarray*}
\left| V (\bbeta,\br) \right|^{2}  
& \leq & \left( n^{ds} \prod_{i=1}^{s} \prod_{\ell=1}^{d} 
	 { \left( \max \left( 1, \left| \alpha_{\ell} \right| \right) \right) }^{2p_{i}(n-1)} 
	 \right)   
	 \times \left( \prod_{j=1}^{d} \prod_{i=0}^{k-1} \left| v_{i} \left( \beta_{j} \right) \right|^{2} 
	 \right) \\
&  =   & n^{ds} M(\alpha)^{2(n-1)(p_{1}+\cdots+p_{s})}   
	 \left( \prod_{j=1}^{d} \prod_{i=0}^{k-1} 
	 \left( \sum_{\ell=i}^{n-1} {\ell \choose i}^{2} \left| \alpha_{j} \right|^{2(\ell-i)}  
	 \right) \right) \\
& \leq & n^{ds} M(\alpha)^{2(n-1)(p_{1}+\cdots+p_{s})}
	 \left( \prod_{j=1}^{d} 
		\max { \left( 1, \left| \alpha_{j} \right| \right) }^{2(n-1)} \right) 
	 { \left( \prod_{i=0}^{k-1}
		\left( \sum_{\ell=i}^{n-1} {\ell \choose i}^{2} \right) 
	   \right) }^{d} \\
&  =   & { \left( n^{s} \prod_{i=0}^{k-1} 
		  \left( \sum_{\ell=i}^{n-1} {\ell \choose i}^{2} \right) 
	   \right) }^{d}        
	 M(\alpha)^{2(n-1)(k+p_{1}+\cdots+p_{s})}.
\end{eqnarray*}     

To complete the proof of the lemma, we now use the fact that 
\begin{displaymath}
\sum_{\ell=i}^{n-1} {\ell \choose i}^{2} 
\leq \sum_{\ell=i}^{n-1} \frac{\ell^{2i}}{(i!)^{2}} 
  =  \frac{1}{(i!)^{2}} \sum_{\ell=i}^{n-1} \ell^{2i}
\leq \frac{1}{(i!)^{2}} \int_{0}^{n} x^{2i}
= \frac{n^{2i+1}}{(2i+1)(i!)^{2}}. 
\end{displaymath}
\end{proof}

\vspace{3.0mm}

{\bf 3. Proof of the Theorem} Proceeding by induction, 
we shall now combine the results we have obtained in 
the previous section to prove our theorem. 

From the work of Boyd and Smyth, the theorem holds for 
$2 \leq d \leq 21$. 

Let $\al$ be an algebraic number of degree $d \geq 22$ 
over $\bbQ$ and assume that the theorem holds for all 
$2 \leq d_{1} < d$. Notice that if $\al$ is not an algebraic 
integer then $a_{d} \geq 2$ and $M(\al) \geq a_{d} \geq 2$. 
Therefore we may assume that $\al$ is an algebraic integer. 

Furthermore, we can also assume that for any prime $p$, 
$\bbQ \left( \al^{p} \right)=\bbQ(\al)$. Otherwise, by Lemme~7.1.1 of
\cite{BP}, the polynomial $f_{p}(X)$ defined in the proof 
of Lemma~\ref{lem:lowbnd} is a power of the minimal 
polynomial of $\al^{p}$. This implies that there exist 
two distinct integers $i$ and $j$ such that $\al_{i}^{p} 
=\al_{j}^{p}$. By Lemme~7.1.2 of \cite{BP}, there exists a 
non-zero algebraic integer, $\beta$, which is not a root of 
unity, of degree less than $d$ with $M(\beta) \leq M(\al)$. 
Since $((\log \log d)/(\log d))^{3}$ is a monotonically 
decreasing function for $d \geq 16$ and since the work of 
Boyd and Smyth shows that 
$M(\beta) > (1/4)((\log \log d)/(\log d))^{3}$ if the 
degree of $\beta$ over $\bbQ$ is less than 22, our inductive 
hypothesis shows that the theorem holds. 

Let us start by comparing the bounds in Lemmas~\ref{lem:lowbnd} 
and \ref{lem:upbnd}; taking the logarithm and dividing both 
sides by $d$, we find that $\log M(\al)$ is at least
\begin{equation}
\label{eq:lowbound}
\frac{\displaystyle 2k \theta(p_{s}) 
+ (k^{2}+s) \left( 3.108 - 8.6d^{-2/3} \right) 
+ \sum_{i=0}^{k-1} \log ((2i+1)(i!)^{2}) - (k^{2}+s) \log (d(k+s))} 
{2(k+s)(k + p_{1} + \cdots + p_{s})}. 
\end{equation}

We will first show that we have 
\begin{displaymath}
0.56 { \left( \frac{\log \log d}{\log d} \right) }^{3} 
< \log M(\alpha), 
\end{displaymath}
for $d \leq 190$, after which we shall show that 
\begin{displaymath}
\frac{1}{4} { \left( \frac{\log \log d}{\log d} \right) }^{3} 
< \log M(\alpha), 
\end{displaymath}
for $d \leq 10000$. 

For $k=7$ and $s=11$, the left-hand side of (\ref{eq:lowbound}) 
is greater than 
\begin{displaymath}
\frac{422.1-516d^{-2/3}-60\log d}{6012}. 
\end{displaymath}

One can compute this quantity for $22 \leq d \leq 94$ to show 
that it is greater than $0.56(\log \log (d)/\log (d))^{3}$. 
Similarly, one can show that the choice $k=8$ and $s=14$ yields 
the desired result for $43 \leq d \leq 190$. Finally, letting 
$k=7$ and $s=17$, we obtain 
\begin{displaymath}
\frac{1}{4} { \left( \frac{\log \log d}{\log d} \right) }^{3} 
< \log M(\alpha), 
\end{displaymath}
for $22 \leq d \leq 10000$. 

For all $d > 10000$, choose 
\begin{displaymath}
k = k_{1} \frac{\log d}{\log \log d} \mbox{ and }
s = s_{1} { \left( \frac{\log d}{\log \log d} \right) }^{2} 
\end{displaymath}
respectively, where the pair $(k_{1},s_{1})$ is contained in some 
region $\cA$ of $\bbR^{2}$. We shall choose this region so that 
(\ref{eq:lowbound}) is at least $(\log \log (d)/ \log (d))^{3}/4$ 
for every pair $(k_{1},s_{1})$ in $\cA$ and also so that there 
is always a pair $(k_{1},s_{1})$ in $\cA$ such that $k$ and 
$s$ are positive integers. For now, we shall assume only that 
$1.2 \leq k_{1} \leq 1.6$ and $1.2 \leq s_{1} \leq 1.51$. 

Since $s_{1} \geq 1.2$ and $d > 10000$, we have $s > 13$, 
and so we may apply Lemma~\ref{lem:primesum}(i) to obtain  
\begin{displaymath}
2k \theta(p_{s}) 
\geq 2k_{1}s_{1} { \left( \frac{\log d}{\log \log d} \right) }^{3} 
     \left( \log s_{1} + 2 \log \log d - 2 \log \log \log d \right) .  
\end{displaymath}

Let us now bound from above the last term in the numerator 
of (\ref{eq:lowbound}). Using the expressions for $k$ and 
$s$ above and the fact that $d > 10000$, we have 
$k + s < 2(\log(d)/ \log \log (d))^{2}$. Thus  
\begin{displaymath}
\log (d(k+s)) < \log d + 2 \log \log d + 0.7 - 2 \log \log \log d. 
\end{displaymath}                               

Using Lemma~\ref{lem:prod} and the above expression for $k$, we have 
\begin{displaymath}
2 \sum_{i=0}^{k-1} \log (i!) 
> { \left( k_{1} \frac{\log d}{\log \log d} \right) }^{2}
  \left( \log \log d - \log \log \log d + \log k_{1} - 1.5 \right) . 
\end{displaymath}

Since $3.108-8.6d^{-2/3} > 3.089$ for $d > 10000$, we can bound 
the last three terms in the numerator of (\ref{eq:lowbound}) 
from below by 
\begin{eqnarray*}
&  & - \frac{(k_{1}^{2}+s_{1}) (\log d)^{3}}{(\log \log d)^{2}} 
     - (k_{1}^{2}+2s_{1}) \frac{(\log d)^{2}}{\log \log d} 
     + (k_{1}^{2} + 2s_{1}) \log \log \log d 
       { \left( \frac{ \log d}{\log \log d} \right) }^{2} \\
&  & + \left( \left( 0.89 + \log k_{1} \right) k_{1}^{2} + 2.39s_{1} \right) 
       { \left( \frac{ \log d}{\log \log d} \right) }^{2}.
\end{eqnarray*}

Notice that to simplify the argument we have ignored the term 
$\sum_{i=0}^{k-1} \log (2i+1)$ which is small for large $d$.  

Applying Lemma~\ref{lem:positive} to this last expression, we 
find that the sum of the last three terms of the numerator of 
(\ref{eq:lowbound}) is greater than 
\begin{displaymath}
- \frac{(k_{1}^{2}+s_{1}) (\log d)^{3}}{(\log \log d)^{2}} 
- \frac{2k_{1}s_{1} \log s_{1} (\log d)^{3}}{(\log \log d)^{2}}. 
\end{displaymath}

Combining this expression with our lower bound for the first term 
(\ref{eq:lowbound}), we find that the numerator of (\ref{eq:lowbound}) 
is greater than 
\begin{displaymath}
\frac{(4k_{1}s_{1}-k_{1}^{2}-s_{1}) (\log d)^{3}}{(\log \log d)^{2}} 
- 4k_{1}s_{1} \log \log \log d 
  { \left( \frac{\log d}{\log \log d} \right) }^{3}. 
\end{displaymath}

Since $\log \log \log d / \log \log d < 0.368$ for $d > 10000$,  
this lower bound yields the still-simpler lower bound 
\begin{equation}
\label{eq:numer}
\frac{(2.528k_{1}s_{1}-k_{1}^{2}-s_{1}) (\log d)^{3}}{(\log \log d)^{2}}.  
\end{equation}

Now let us obtain an upper bound for the denominator of 
(\ref{eq:lowbound}). Using the expressions for $k$ and $s$ 
along with the bound in Lemma~\ref{lem:primesum}(iii), which 
is possible since we saw above that $s > 9$, we find that this 
denominator is at most 
\begin{equation}
\label{eq:denomprod}
\frac{2.256s_{1}^{3}(\log d)^{6}}{(\log \log d)^{5}}
\left( 1 + \frac{k_{1} \log \log d}{s_{1} \log d} \right) 
\left( \frac{k_{1} (\log \log d)^{2}}{1.128s_{1}^{2}(\log d)^{3}} 
	+ 1 + \frac{\log s_{1} - 2\log \log \log d}
		   {2 \log \log d} 
\right) .
\end{equation}

We shall show that the product of the last two terms is less than 
1. To demonstrate this we shall prove that the function 
\begin{eqnarray*}
g(k_{1},s_{1},d) 
& = & \frac{2 \log \log \log d - \log s_{1}}{2 \log \log d} 
      - \frac{k_{1} (\log \log d)^{2}}{1.128s_{1}^{2}(\log d)^{3}} 
      - \frac{k_{1} \log \log d}{s_{1} \log d} \\ 
& = & \left( 282s_{1}^{2} \log^{3}d \log \log \log d 
	     - 141 s_{1}^{2} \log s_{1} \log^{3} d  \right. \\
&   & \left. - 282 k_{1} s_{1} \log^{2} d ( \log \log d)^{2} 
	     - 250 k_{1} ( \log \log d)^{3} \right)  
	   / \left( 282 s_{1}^{2} \log^{3} d \log \log d \right)
\end{eqnarray*}
is positive for $d > 10000, 1.2 \leq s_{1} \leq 1.51$ and 
$1.2 \leq k_{1} \leq 1.05s_{1}$. This will imply that the product 
of the last two terms in (\ref{eq:denomprod}) can be written as 
$(1+a)(1-b)$ where $a$ and $b$ are functions of $k_{1},s_{1}$ 
and $d$ which satisfy $b > a$, from which our claim follows. 

We need only examine the numerator of the expression for $g$. 
Since $k_{1}$ only occurs in terms being subtracted, we may 
replace it by $1.05s_{1}$ in light of our constraint on $k_{1}$. 
Therefore we need only show that 
\begin{displaymath}
2820s_{1} \log^{3}d \log \log \log d - 1410 s_{1} \log s_{1} \log^{3} d
- 2961 s_{1} \log^{2} d ( \log \log d)^{2} - 2625 ( \log \log d)^{3} 
\end{displaymath}
is positive for $d > 10000$ and $1 \leq s_{1} \leq 1.51$, 
upon dividing by $s_{1}$. 

We can rewrite this expression as       
\begin{displaymath}
s_{1} \log^{3} d 
\left( 2820 \log \log \log d - 1410 \log s_{1} 
       - \frac{2961 ( \log \log d)^{2}}{\log d} \right) 
- 2625 ( \log \log d)^{3}.   
\end{displaymath}

Now since $(\log \log d)^{2}/ \log d < 0.54$ for $d > 10000$ 
and since $s_{1} \leq 1.51$, this expression is greater than 
\begin{displaymath}
69s_{1} \log^{3} d - 2625 ( \log \log d)^{3}.   
\end{displaymath}

Since $\log \log d/ \log d < 0.25$ for $d > 10000$, this last 
expression is greater than $41 \log^{3} d$ for $s_{1} \geq 1.2$ 
and $d > 10000$. Hence $g(d,k_{1},s_{1})$ is positive in the 
desired region. 

Therefore the product in question is less than 1 which implies that 
\begin{equation}
\label{eq:denom}
2(k+s)(k+p_{1}+\cdots+p_{s}) 
\leq 2.256s_{1}^{3} \frac{(\log d)^{6}}{(\log \log d)^{5}}
\end{equation}

Combining \eqref{eq:numer} with the upper bound in \eqref{eq:denom}, 
we will show that 
\begin{displaymath}
\frac{2.528k_{1}s_{1}-k_{1}^{2}-s_{1}}
     {2.256s_{1}^{3}} 
> \frac{1}{4}      
\end{displaymath}
for all $1.26 \leq k_{1} \leq 1.51$ and $k_{1}-0.06 \leq s_{1} \leq k_{1}$.  
Notice that this choice of $k_{1}$ and $s_{1}$ satisfies the 
conditions set down above. Also, since $d > 10000$, there 
exist $k_{1}$ and $s_{1}$ in this range such that $k$ and $s$ 
are integers. 

Let us now turn to proving this inequality. We shall use ideas 
from multi-variable calculus. Let 
\begin{displaymath}
f(k_{1},s_{1}) = \frac{2.528k_{1}s_{1}-k_{1}^{2}-s_{1}}
		      {2.256s_{1}^{3}}.  
\end{displaymath}

Then 
\begin{displaymath}
\frac{\partial f}{\partial k_{1}} 
= \frac{158s_{1}-125k_{1}}{141s_{1}^{3}} \mbox{ and }
\frac{\partial f}{\partial s_{1}} 
= \frac{375k_{1}^{2}-632k_{1}s_{1}+250s_{1}}{282s_{1}^{4}}. 
\end{displaymath}

From these expressions, we find that both derivatives are 
zero only at $k_{1}=125/79=1.58 \ldots$ and $s_{1}=15625/12482
=1.25\ldots$ --- but these values do not fall inside our 
range and hence there are no local minima for such $k_{1}$ 
and $s_{1}$. Thus we need only look along the boundary to 
determine the minimum of $f$ for $(k_{1},s_{1})$ in this 
region. Along the edge formed by the line $k_{1}=1.26$, we 
have $f(1.26,s_{1})=(27316s_{1}-19845)/(28200s_{1}^{3})$. 
For $1.2 \leq s_{1} \leq 1.26$, this function attains its 
minimum, which is greater than $0.2583$, at $s_{1}=1.26$. 
Along the edge formed by the line $k_{1}=1.51$, we have 
$f(1.51,s_{1})=(140864s_{1}-114005)/(112800s_{1}^{3})$. 
For $1.45 \leq s_{1} \leq 1.51$, the minimum of this function 
occurs at $s_{1}=1.51$ and is greater than $0.2541$. Along the 
boundary formed by the line $k_{1}=s_{1}+0.06$, $f(k_{1},s_{1})
=(19100s_{1}^{2}-12104s_{1}-45)/(28200s_{1}^{3}) \geq 0.2624$ 
--- this value occurs for $s=1.45$.  Finally, along the boundary 
formed by the line $k_{1}=s_{1}$, $f(k_{1},s_{1})=(191s_{1}-125)
/(282s_{1}^{2}) \geq 0.2541$ --- this value occurs for $s=1.51$.  

This completes the proof of the theorem. 

\setlength{\baselineskip}{5 mm}

\noindent
DEPARTMENT OF MATHEMATICS \\ 
CITY UNIVERSITY \\ 
NORTHAMPTON SQUARE \\ 
LONDON EC1V 0HB \\
UK

\end{document}